\documentclass[11pt]{article} 
\usepackage{exscale,relsize}
\usepackage{amsmath}
\usepackage{amsfonts}
\usepackage{amssymb}
\usepackage{calc}
\usepackage{theorem}
\usepackage{pifont}      
\usepackage{graphicx}
\newcommand{\emp}{\ensuremath{\varnothing}}

\newcommand{\scal}[2]{\langle{{#1},{#2}}\rangle}

\newcommand{\exi}{\ensuremath{\exists\,}}

\newcommand{\RR}{\ensuremath{\mathbb R}}

\newcommand{\RP}{\ensuremath{\left[0,+\infty\right[}}

\newcommand{\RX}{\ensuremath{\,\left]-\infty,+\infty\right]}}

\newcommand{\NN}{\ensuremath{\mathbb N}}
\newcommand{\nnn}{\ensuremath{{n \in \NN}}}
\newcommand{\thalb}{\ensuremath{\tfrac{1}{2}}}

\newcommand{\menge}[2]{\big\{{#1} \mid {#2}\big\}}

\newcommand{\To}{\ensuremath{\rightrightarrows}}

\newcommand{\dom}{\ensuremath{\operatorname{dom}}}

\newcommand{\gra}{\ensuremath{\operatorname{gra}}}

\newcommand{\ran}{\ensuremath{\operatorname{ran}}}

\newcommand{\cldom}{\ensuremath{\overline{\operatorname{dom}}\,}}

\newcommand{\Id}{\ensuremath{\operatorname{Id}}}

\newcommand{\minf}{\ensuremath{-\infty}}
\newcommand{\pinf}{\ensuremath{+\infty}}

\renewcommand{\phi}{\ensuremath{\varphi}}

\newcommand{\halb}{\ensuremath{\tfrac{1}{2}}}

\newtheorem{theorem}{Theorem}[section]

\newtheorem{fact}[theorem]{Fact}
\newtheorem{corollary}[theorem]{Corollary}
\newtheorem{proposition}[theorem]{Proposition}
\newtheorem{definition}[theorem]{Definition}

\theoremstyle{plain}{\theorembodyfont{\rmfamily}
}
\theoremstyle{plain}{\theorembodyfont{\rmfamily}
}
\theoremstyle{plain}{\theorembodyfont{\rmfamily}
}
\theoremstyle{plain}{\theorembodyfont{\rmfamily}
\newtheorem{example}[theorem]{Example}}
\theoremstyle{plain}{\theorembodyfont{\rmfamily}
\newtheorem{remark}[theorem]{Remark}}
\theoremstyle{plain}{\theorembodyfont{\rmfamily}
}

\begin{document}

\title{{\sffamily Autoconjugate representers for 
linear monotone operators}}

\author{
Heinz H.\ Bauschke\thanks{Mathematics, Irving K.\ Barber School,
UBC Okanagan, Kelowna, British Columbia V1V 1V7, Canada. E-mail:
\texttt{heinz.bauschke@ubc.ca}.},~  Xianfu
Wang\thanks{Mathematics, Irving K.\ Barber School, UBC Okanagan,
Kelowna, British Columbia V1V 1V7, Canada. E-mail:
\texttt{shawn.wang@ubc.ca}.},~ and Liangjin\
Yao\thanks{Mathematics, Irving K.\ Barber School, UBC Okanagan,
Kelowna, British Columbia V1V 1V7, Canada.
E-mail:  \texttt{ljinyao@interchange.ubc.ca}.}. }

\date{February 10, 2008}

\maketitle


\begin{abstract} \noindent
Monotone operators are of central importance in modern optimization and
nonlinear analysis. Their study has been revolutionized lately, due to the
systematic use of the Fitzpatrick function. Pioneered by Penot and Svaiter, 
a topic of recent
interest has been the representation of maximal monotone operators
by so-called autoconjugate functions. 
Two explicit constructions were proposed, the first by Penot and
Z\u{a}linescu in 2005, and another by Bauschke and Wang in 2007.
The former requires a mild constraint qualification while the latter is
based on the proximal average. 

We show that these two autoconjugate representers must coincide for
continuous linear monotone operators on reflexive spaces. 
The continuity and
the linearity assumption are both essential as examples of discontinuous
linear operators and of subdifferential operators illustrate. 
Furthermore, we also construct an infinite family of 
autoconjugate representers 
for the identity operator on the real line.
\end{abstract}



\section{Introduction}

Throughout this paper, we assume that
\begin{equation*}
\text{$X$ is a real reflexive Banach space, with continuous dual space
$X^*$, and pairing $\scal{\cdot}{\cdot}$.}
\end{equation*}
The norm of $X$ is denoted by $\|\cdot\|$, and the norm in the dual space
$X^*$
by $\|\cdot\|_*$.

Let $A\colon X\To X^*$ be a set-valued operator, with
\emph{graph} 
$\gra A = \menge{(x,x^*)\in X\times X^*}{x^*\in Ax}$,
with \emph{inverse operator} $A^{-1}\colon X^*\To X$ given
by $\gra A^{-1} = \menge{(x^*,x)\in X^*\times X}{x^*\in Ax}$,
with \emph{domain} $\dom A = \menge{x\in X}{Ax\neq\emp}$, and
with \emph{range} $\ran A = A(X)$. 
Recall that $A$ 
is \emph{monotone} if
\begin{equation}
\big(\forall (x,x^*)\in \gra A\big)\big(\forall (y,y^*)\in\gra A\big)
\quad \scal{x-y}{x^*-y^*}\geq 0. 
\end{equation}
A monotone operator $A$ is \emph{maximal
monotone} if no proper enlargement (in the sense of graph inclusion) 
of $A$ is monotone.
Monotone operators are ubiquitous in Optimization and Analysis 
(see, e.g.,  \cite{BuIu,Phelps,RockWets,Si,SiNewBook,Zalinescu}) 
since they contain the key classes of subdifferential
operators and of positive linear operators. 

In \cite{FIT}, Fitzpatrick introduced the following 
tool in the study of monotone operators. 

\begin{definition}
Let $A\colon X\rightrightarrows X^*.$ The
\emph{Fitzpatrick function} of $A$ is
\begin{equation}
F_A\colon  (x,x^*)\mapsto \sup_{(y,y^*)\in
\gra A}\langle x, y^*\rangle+ \langle y,x^*\rangle-\langle y,y^*\rangle.
\end{equation}
\end{definition}
Monotone Operator Theory has been revolutionized 
through the systematic use of the Fitzpatrick function; new results
have been obtained and previously known result
have been reproved in a simpler fashion --- see, e.g., 
\cite{BBBRW, BPAMS, BBW, BLW, BM, BaW, BC, Borwein06-, Borwein06,
BorweinZhu, BCW, BGW, BuFi, BuIu06, GLR, Gh, Msv, MT, Penot, Penot2,
PenotZ, RS, Si06, Si06b, Si07, SiZ-, SiZ, Sv, Voi06, Zali05}. 
Before listing some of the key properties 
of the Fitzpatrick function, we introduce a convenient notation
utilized by Penot \cite{Penot2}:
If $F\colon X\times X^*\to\RX$, set 
\begin{equation}
F^\intercal\colon X^*\times X\colon (x^*,x)\mapsto F(x,x^*), 
\end{equation}
and similarly for a function defined on $X^*\times X$. 
We now define an associated operator $X\To X^*$ by requiring that 
for $(x,x^*)\in X\times X^*$, 
\begin{equation}
x^* \in  G(F)x \quad \Leftrightarrow \quad
F(x,x^*) = \scal{x}{x^*};
\end{equation}
we also say that $F$ is a \emph{representer} for $G(F)$.

\begin{fact} \label{f:Fitz} {\rm (See \cite{FIT}.)} 
Let $A\colon X\To X^*$ be maximal monotone. Then the following
hold.
\begin{enumerate}
\item \label{f:Fitz:F1} $F_A$ is proper, lower semicontinuous, and convex. 
\item \label{f:Fitz:ndif:1} $F_{A^{-1}} = F_A^\intercal$.  
\item \label{f:Fitz:gep} $F_A \geq \scal{\cdot}{\cdot}$.
\item \label{f:Fitz:rep} $A = G(F_A)$. 
\end{enumerate}
\end{fact}
Item~\ref{f:Fitz:rep} of Fact~\ref{f:Fitz} states the key property
that the Fitzpatrick function $F_A$ is a representer for the maximal monotone
operator $A$. 
It turns out that there are even more structured 
representers for $A$ available: recall that 
$F\colon X\times X^*\to \RX$ is \emph{autoconjugate}, if
\begin{equation}
F^* = F^\intercal.
\end{equation}

Autoconjugate representers are readily available for two
important classes of maximal monotone operators.

\begin{example}[subdifferential operator] \label{ex:subdiff} 
Let $f\colon X\to\RX$ be proper, lower semicontinuous, and convex.
Then the separable sum of $f$ and the Fenchel conjugate $f^*$, i.e., 
\begin{equation}
f\oplus f^* \colon X\times X^* \to \RX \colon (x,x^*) \mapsto f(x)+f^*(x^*),
\end{equation}
is an autoconjugate representer for the subdifferential operator $\partial f$. 
\end{example}

\begin{example}[antisymmetric operator] \label{ex:antisymm}
Let $A\colon X\to X^*$ be continuous, linear, and \emph{antisymmetric}, 
i.e., $A^* = -A$.
Then the indicator function of the graph of $A$, i.e., 
\begin{equation}
\iota_{\gra A} \colon X\times X^*\to\RX\colon
(x,x^*)\mapsto
\begin{cases}
0, &\text{if $x^*=Ax$;}\\
\pinf, &\text{otherwise} 
\end{cases}
\end{equation}
is an autoconjugate representer for $A$. 
\end{example}

We now list some very pleasing and well known properties of 
autoconjugate functions. 

\begin{fact}[Penot-Simons-Z\u{a}linescu] 
{\rm (See  \cite{Penot,Penot2,PenotZ,SiZ}.)}
\label{f:auto} 
Let 
$F\colon X\times X^*\to \RX$ be autoconjugate.
Then the following hold.
\begin{enumerate}
\item $F$ is proper, lower semicontinuous, and convex. \label{f:auto:Unl:4}
\item $F\geq \scal{\cdot}{\cdot}$. \label{Aum:1}
\item $G(F)$ is maximal monotone. \label{f:auto:GF}
\item If $\widetilde{F}\colon X\times X^*\to\RX$ is autoconjugate
and $\widetilde{F}\leq F$, then $\widetilde{F}=F$. \label{Au:4}
\end{enumerate}
\end{fact}

Unfortunately, the Fitzpatrick function $F_A$ is 
usually \emph{not} an autoconjugate
representer for $A$. In view of Fact~\ref{f:auto}\ref{f:auto:GF},
it was tempting to ask whether every general maximal monotone operator
possesses an autoconjugate representer.
Nonconstructive existence proofs were presented 
by Penot \cite{Penot,Penot2} and by
Svaiter \cite{Sv} in 2003 (see also Ghoussoub's preprint \cite{Gh}). 
The first actual construction of an autoconjugate representer for a
maximal monotone operator satisfying a mild constraint qualification was
provided by Penot and Z\u{a}linescu in 2005. 

\begin{fact}[Penot-Z\u{a}linescu] {\rm (See \cite{PenotZ}.)}
\label{GFF:2}
Let $A\colon X\To X^*$ be maximal monotone. Suppose that
the affine hull of $\dom A$ is closed. Then 
\begin{align} 
\mathcal{A}_A\colon X\times X^* &\to \RX\nonumber\\
(x,x^*) &\mapsto \inf_{y^*\in X^*}
\tfrac{1}{2}F_A(x,x^*+y^*)+\tfrac{1}{2}F_A^{*\intercal}(x,x^*-y^*)\label{e:18}
\end{align}
is an autoconjugate representer for $A$. 
\end{fact}

Another autoconjugate representer was very recently proposed in \cite{BaW}. 
While this proximal-averaged based construction is 
more involved \cite{BR,BLM,BMA}, 
it has the advantage of not imposing a constraint qualification.

\begin{fact} {\rm (See \cite{BaW}.)}
\label{GFF:1}
Let $A\colon X\To X^*$ be maximal monotone. Then
\begin{align}
\mathcal{B}_A \colon X\times X^*&\to \RX\nonumber\\
 (x,x^*)&\mapsto 
\inf_{(y,y^*)\in X\times X^* }
\tfrac{1}{2}F_A(x+y,x^*+y^*)+
\tfrac{1}{2}F^{*\intercal}_A(x-y,x^*-y^*)+\tfrac{1}{2}\|y\|^2
+\tfrac{1}{2}\|y^*\|_*^2 \label{e:GFF:1}
\end{align}
is an autoconjugate representer for $A$.
\end{fact}

It is natural to ask ``How do the autoconjugate representers
$\mathcal{A}_A$ and $\mathcal{B}_A$ compare?''
We provide two answers to this question:
First, we show that if $A\colon X\to X^*$ is continuous, linear, and monotone,
then $\mathcal{A}_A$ and $\mathcal{B}_A$ coincide; furthermore,
we provide a formula for this autoconjugate representer 
which agrees with a third autoconjugate representer $\mathcal{C}_A$ 
that is contained in the work by Ghoussoub (Theorem~\ref{t:main}). 
Secondly, for \emph{nonlinear} monotone subdifferential operators,
the two 
autoconjugate representers may be different (Theorem~\ref{t:different}). 

The first answer raises the question on whether autoconjugate representers
for continuous linear monotone operators are unique. We answer this
question in the negative by providing a family of autoconjugate
representers for the identity operator $\Id$ (Theorem~\ref{t:rrsp}). 
However, we show that the autoconjugate representers $\mathcal{A}_A$ and
$\mathcal{B}_A$ in this setting are
characterized by a pleasing symmetry property (Theorem~\ref{t:hoe}).

We conclude by discussing  \emph{discontinuous} linear monotone operators.
It turns out that $\mathcal{A}_A$ may fail to be autoconjugate
(Example~\ref{ex:potato}), which
underlines not only the continuity assumption in Theorem~\ref{t:main} 
but also the importance of the constraint qualification in
Fact~\ref{GFF:2}. 

The remainder of this paper is organized as follows.
Section~\ref{s:auxi} contains some results on quadratic functions and
another autoconjugate representer that will be used in later sections.
In Section~\ref{s:coin}, we show that $\mathcal{A}_A$ and $\mathcal{B}_A$
coincide and provide a simple formula for it.
Uniqueness of autoconjugate representations are discussed in
Section~\ref{s:uniq}, and a characterization in the symmetric case is also
presented. In stark contrast, and as shown in Section~\ref{s:-ln},
$\mathcal{A}_A$ and $\mathcal{B}_A$ may be different for (nonlinear)
subdifferential operators. The final Section~\ref{s:disc} reveals similar
difference for discontinuous linear operators. 

Notation utilized is standard as in Convex Analysis and Monotone
Operator Theory; see, e.g., \cite{Rocky,RockWets,Zalinescu}. 
Thus, for a proper convex function $f\colon X\to\RX$, we write
$f^* \colon x^*\mapsto \sup_{x\in X} \scal{x}{x^*}-f(x)$,~
$\partial f\colon X\To X^*\colon x\mapsto \menge{x^*\in X^*}{(\forall y\in
X)\; \scal{y-x}{x^*} + f(x)\leq f(y)}$,~  and 
$\dom f = \menge{x\in X}{f(x)<\pinf}$, 
for the \emph{Fenchel conjugate}, \emph{subdifferential operator}, 
and \emph{domain} of $f$, respectively. 
The strictly positive integers are $\NN = \{1,2,\ldots\}$.

\section{Auxiliary Results}

\label{s:auxi}

The following result in a consequence of results and proof techniques
introduced by Penot, Simons, and Z\u{a}linescu \cite{PenotZ, SiZ}.
It also extends \cite[Lemma~2.2]{Gh}.

\begin{proposition}\label{sumoperation}
Let $F_{1}$ and $F_2$ be autoconjugate functions on $X\times X^*$ 
representing maximal monotone operators $A_1$ and $A_2$,
respectively. Suppose that 
\begin{equation}\label{CQ}
\bigcup_{\lambda>0}\lambda \big(P_{X}\dom F_{1}-P_{X}\dom F_{2}\big) 
\quad \text{is a closed subspace of $X$,}
\end{equation}
where $P_{X}\colon X\times X^*\to X\colon (x,x^*)\mapsto x$, 
and set 
\begin{equation} \label{e:schnee}
F\colon X\times X^*\to\RX \colon 
(x,x^*)\mapsto \inf_{y^*\in X^*} F_{1}(x,y^*)+F_{2}(x,x^*-y^*).
\end{equation}
Then $F$ is an autoconjugate representer for $A_1+A_2$, and
the infimum in \eqref{e:schnee} is attained. 
\end{proposition}
\begin{proof}
Let $(x,x^*)\in X\times X^*$. 
Using Simons and Z\u{a}linescu's \cite[Theorem~4.2]{SiZ} and
the assumption that each $F_i$ is autoconjugate, we obtain
\begin{align}
F^*(x^*,x) &= 
\min_{x_1^*+x_2^*=x^*}\; F^*_{1}(x^*_{1},x)+F_{2}^*(x^*_{2},x)\notag\\
&=\min_{x_1^*+x_2^*=x^*}\; F_{1}(x,x^*_{1})+F_{2}(x,x^*_{2})\notag\\
&= F(x,x^*).
\end{align}
Thus, $F$ is autoconjugate and the infimum in \eqref{e:schnee} is attained. 

It remains to show that $G(F) = G(F_1)+G(F_2)$. 
Since autoconjugates are greater than or equal to 
$\scal{\cdot}{\cdot}$ (see Fact~\ref{f:auto}\ref{Aum:1}),
the above implies the equivalences
\begin{align}
x^*\in G(F)x &\Leftrightarrow F(x,x^*)=\scal{x}{x^*}\notag\\
&\Leftrightarrow (\exi y^*\in X^*)\;\; F_1(x,y^*) + F_2(x,x^*-y^*) =
\scal{x}{y^*} + \scal{x}{x^*-y^*}\notag\\
&\Leftrightarrow (\exi y^*\in X^*)\;\;
F_1(x,y^*)=\scal{x}{y^*}\;\;\text{and}\;\; F_2(x,x^*-y^*) =
\scal{x}{x^*-y^*}\notag\\
&\Leftrightarrow (\exi y^*\in X^*)\;\;
y^*\in G(F_1)(x)\;\;\text{and}\;\; x^*-y^*\in G(F_2)(x)\notag\\
&\Leftrightarrow (\exi y^*\in X^*)\;\;
y^*\in A_1x\;\;\text{and}\;\; x^*-y^*\in A_2x\notag\\
&\Leftrightarrow  x^*\in(A_1+A_2)x.
\end{align}
Therefore, $G(F)=A_1+A_2$, i.e., $F$ is a representer for $A_1+A_2$. 
\end{proof}

Suppose that
\begin{equation}
A \colon X\to X^* \quad\text{is linear and continuous.}
\end{equation}
Then $A$ is \emph{symmetric} (resp.\ \emph{antisymmetric}) if
$A^*=A$ (resp.\ $A^*=-A$). We denote the \emph{symmetric part}
and the  \emph{antisymmetric part}
of $A$ by
\begin{equation}
A_+ = \thalb A + \thalb A^* \quad\text{and}\quad
A_{\mathlarger{\circ}} = \thalb A - \thalb A^*, 
\end{equation}
respectively. Throughout, we shall work with the quadratic function
\begin{equation} \label{e:diequad}
q_A \colon X\to\RR\colon x\mapsto \thalb \scal{x}{Ax}, 
\end{equation}
and we will use the well known facts (see, e.g., \cite{PheSim}) that
$q_A = q_{A_+}$, that 
\begin{equation} \label{e:gradq}
\nabla q_A = A_+,
\end{equation}
and that 
$A$ is maximal monotone if and only if $q_A$ is convex. 
The next result provides a formula for $q_A^*$ that will be useful later. 

\begin{proposition}\label{better} 
Let $A\colon X\to X^*$ be continuous, linear, symmetric, and monotone. 
Then 
\begin{equation} \label{e:better}
\big(\forall (x,x^*)\in X\times X^*\big)\quad
q_{A}^*(x^*+Ax)=q_{A}(x)+\scal{x}{x^*} +q_{A}^*(x^*) 
\end{equation}
and
\begin{equation} \label{Alf:1}
q_A^* \circ A = q_A. 
\end{equation}

\end{proposition}
\begin{proof}
Let $(x,x^*)\in X\times X^*$. Then 
\begin{align}
q^*_A(x^*+Ax)&=\sup_{y}\; \scal{y}{x^*+Ax}-q_A(y)=
\sup_{y}\; \scal{y}{x^*}-q_A(y)+\scal{y}{Ax}\notag\\
&= q_A(x) + \sup_{y}\; \scal{y}{x^*}-q_A(y)+
\scal{y}{Ax}-q_A(x)
= q_A(x)+ \sup_{y}\;\scal{y}{x^*}-q_A(y-x)\notag\\
&=q_A(x)+ \scal{x}{x^*} + \sup_{y}\;\scal{y-x}{x^*}-q_A(y-x)
=q_A(x) + \scal{x}{x^*} + q^*_A(x^*),
\end{align}
which verifies \eqref{e:better}. To see \eqref{Alf:1}, set 
$x^*=0$ in \eqref{e:better}. 
\end{proof}

\begin{proposition}
Let $A\colon X\to X^*$ be continuous, linear, and monotone, and 
let $(x,x^*)\in X\times X^*$.
Then 
\begin{equation} \label{e:es:4}
F_A(x,x^*)=2
q_{A_+}^*(\tfrac{1}{2}x^*+\tfrac{1}{2}A^*x)=\tfrac{1}{2}q_{A_{+}}^*(x^*+A^*x).
\end{equation}
and 
\begin{equation}
\label{e:Faz:1}
F_A^*(x^*,x) = \iota_{\gra A}(x,x^*) + \scal{x}{Ax}.
\end{equation}
\end{proposition}
\begin{proof}
As in the proof of \cite[Theorem~2.3(i)]{BBW}, we have
\begin{align}
F_A(x,x^*) &= \sup_{y\in X} \scal{x}{Ay} +\scal{y}{x^*} -
\scal{y}{Ay}\notag\\
&= 2\sup_{y\in X} \scal{y}{\thalb x^* + \thalb A^*x} - q_{A_+}(y)\notag\\
&= 2q_{A_+}^*(\thalb x^*+\thalb A^*x)\notag\\
&= \thalb q_{A_+}^*(x^* + A^*x).
\end{align}
This verifies \eqref{e:es:4}. Furthermore, \eqref{e:Faz:1} follows from
$F_A^*(x^*,x) = (\iota_{\gra A}+\scal{\cdot}{\cdot})^{*\intercal *}(x^*,x)
=\iota_{\gra A}(x,x^*) + \scal{x}{Ax}$. 
\end{proof}

\begin{proposition}\label{PAD}
Let $F_1\colon X\times X^*\to\RX$ be autoconjugate,
and let 
$A_2\colon X\to X^*$ be continuous, linear, and antisymmetric. 
Then the function 
\begin{equation}
(x,x^*)\mapsto F_1(x,x^*-A_2x)
\end{equation}
is an autoconjugate representer for $G(F_1)+A_2$. 
\end{proposition}
\begin{proof} 
Set $F_2 = \iota_{\gra A_2}$.
By Example~\ref{ex:antisymm}, $F_2$ is an autoconjugate representer for
$A_2$.
Let $F$ be as in Proposition~\ref{sumoperation}. Then for every $(x,x^*)\in
X\times X^*$, we have
\begin{align}
F(x,x^*) &= \inf_{z^*\in X^*} F_1(x,x^*-z^*) + F_2(x,z^*)\notag\\
&= \inf_{z^*\in X^*} F_1(x,x^*-z^*) + \iota_{\gra A_2}(x,z^*)\notag\\
&= F_1(x,x^*-A_2x).
\end{align}
Thus, Proposition~\ref{sumoperation} yields that $F$ represents
$G(F_1)+A_2$. 
\end{proof}

\begin{example}[Ghoussoub]  \label{ex:Gh} 
{\rm (See also \cite[Section~1]{Gh}.)} 
Let $f\colon X\to\RX$ be proper, lower semicontinuous, and convex,
and let $A$ be antisymmetric. Then 
the function
\begin{equation}
(x,x^*)\mapsto f(x) + f^*(x^*-Ax)
\end{equation}
is an autoconjugate representer for $\partial f + A$. 
\end{example}
\begin{proof}
By Example~\ref{ex:subdiff}, 
$f\oplus f^*$ is an autoconjugate representer for
$\partial f$. The result thus follows from Proposition~\ref{PAD}. 
\end{proof}

\begin{corollary}\label{c:G}
Let $A\colon X \to X^*$ be continuous, linear, and
monotone. Then 
\begin{align} 
\mathcal{C}_A\colon X\times X^* &\to\RX\nonumber\\
(x,x^*) &\mapsto q_{A_+}(x) + q^*_{A_+}(x^*-A_{\mathlarger \circ}x)
\label{e:G}
\end{align}
is an autoconjugate representer for $A$. 
In particular, if $A$ is symmetric, then
\begin{equation} \label{e:symmG}
\mathcal{C}_A = q_A \oplus q_A^*.
\end{equation}
\end{corollary}
\begin{proof}
This follows from \eqref{e:gradq} and 
Example~\ref{ex:Gh} (when applied to the function $q_A=q_{A_+}$ and 
to the antisymmetric operator $A_{\mathlarger{\circ}}$). 
\end{proof}

We now show that the Ghoussoub representers are closed under the partial
infimal convolution operation of Proposition~\ref{sumoperation}. 

\begin{proposition} 
Let $A$ and $B$ be continuous, linear, and monotone on $X$.
Then the function
\begin{equation} \label{e:bigwhite}
F\colon X\times X^*\to\RX\colon
(x,x^*)\mapsto \min_{y^*\in X^*} \mathcal{C}_{A}(x,x^*-y^*) +
\mathcal{C}_{B}(x,y^*)
\end{equation}
coincides with the autoconjugate representer $\mathcal{C}_{A+B}$ for
$A+B$.
\end{proposition}

\begin{proof}
In view of Proposition~\ref{sumoperation},
we only need to show that $F=\mathcal{C}_{A+B}$.
Let $(x,x^*)\in X\times X^*$.
Using \eqref{e:bigwhite} and  Corollary~\ref{c:G},
we obtain
\begin{align}
F(x,x^*) 
&= \min_{y^*\in X^*} q_{A_+}(x) + q_{A_+}^*(x^*-y^*-A_{\mathlarger\circ}x) +
q_{B_+}(x) + q_{B_+}^*(y^*-B_{\mathlarger\circ}x)\nonumber\\
&=q_{A_+}(x)+q_{B_+}(x) + (q_{A_+}^*\Box
q_{B_+}^*)(x^*-A_{\mathlarger\circ}x-B_{\mathlarger\circ}x)\nonumber\\
&=q_{A_++B_+}(x) +
\big(q_{A_+}+q_{B_+}\big)^*(x^*-A_{\mathlarger\circ}x-B_{\mathlarger\circ}x)\nonumber\\
&= q_{(A+B)_+}(x) +
q_{(A+B)_+}^*\big(x^*-(A+B)_{\mathlarger\circ}x\big)\nonumber\\
&= \mathcal{C}_{A+B}(x,x^*),
\end{align}
as required.
\end{proof}

\section{Coincidence} 

\label{s:coin}

We are now ready for one of our main results. 

\begin{theorem}[coincidence] \label{t:main}
Let $A\colon X\to X^*$ be continuous, linear, and monotone.
Then all three autoconjugate representers $\mathcal{A}_A$, $\mathcal{B}_A$,
$\mathcal{C}_A$ for $A$ coincide with the function
\begin{equation} \label{e:29}
(x,x^*)\mapsto 
\scal{x}{x^*} + q_{A_+}^*(x^*-Ax).
\end{equation}
\end{theorem}
\begin{proof}
The proof proceeds by proving a succession of claims. Let $(x,x^*)\in
X\times X^*$. 

\emph{Claim~1:} $\mathcal{A}_A = \mathcal{C}_A$.\\
Using \eqref{e:18}, \eqref{e:Faz:1}, and \eqref{e:G}, we obtain
\begin{align}
\mathcal{A}_A(x,x^*)&= 
\inf_{y^*\in X^*}\halb F_A\big(x,x^*+y^*)\big)+ \thalb F_A^*(x^*-y^*,x)\nonumber\\
&=\inf_{y^*\in X^*} \halb F_A\big(x,x^*+y^*)\big)+
\iota_{\gra A}(x,x^*-y^*) +q_{A}(x) \nonumber\\
&=\thalb F_A(x,2x^*-Ax) +q_{A_+}(x)\nonumber\\
&= q_{A_+}^*\big(x^*-\thalb Ax+ \thalb A^*x\big) + q_{A_+}(x)
\label{mainequation}\\
&= q_{A_+}^*\big(x^*-A_{\mathlarger \circ}x\big) + q_{A_+}(x)\nonumber\\
&= \mathcal{C}_A(x,x^*).
\end{align}
This verifies Claim~1.

\emph{Claim~2:} $\mathcal{A}_A$ coincides with the function of \eqref{e:29}.\\
In view of \eqref{mainequation} and Proposition~\ref{better}, we see that 
\begin{align}
\mathcal{A}_A(x,x^*)&= 
q_{A_+}^*\big(x^*-\thalb Ax+ \thalb A^*x\big) + q_{A_+}(x)\nonumber\\
&=q_{A_+}^*(x^*-Ax+A_+x)+q_{A_+}(x)\nonumber\\
&=2 q_{A_+}(x) +\scal{x}{x^*-Ax} + q_{A_+}^*(x^*-Ax)\nonumber\\
&= \scal{x}{x^*} + q_{A_+}^*(x^*-Ax),
\end{align} 
which establishes Claim~2.

\emph{Claim~3:} $\mathcal{A}_A = \mathcal{B}_A$.\\
Using \eqref{e:GFF:1}, \eqref{e:Faz:1}, \eqref{e:es:4}, 
Proposition~\ref{better}, and Claim~2, we have 
{\allowdisplaybreaks
\begin{align} \mathcal{B}_A(x,x^*) 
&= \inf_{(y,y^*)\in X\times X^*}
\tfrac{1}{2}F_A(x+y,x^*+y^*)+ \tfrac{1}{2}F^{*\intercal}_A(x-y,x^*-y^*)+\tfrac{1}{2}\|y\|^2 +\tfrac{1}{2}\|y^*\|_*^2\nonumber \\
&= \inf_{(y,y^*)\in X\times X^* }
\tfrac{1}{2}F_A(x+y,x^*+y^*)+\iota_{\gra A}(x-y,x^*-y^*)  \nonumber\\
&\qquad\qquad\qquad +
\tfrac{1}{2}\scal{x-y}{A(x-y)}+\tfrac{1}{2}\|y\|^2 +\tfrac{1}{2}\|y^*\|_*^2
\nonumber\\
&=\inf_{y\in X }
\tfrac{1}{2}F_A(x+y,2x^*-A(x-y))+q_A(x-y) +\tfrac{1}{2}\|y\|^2+\tfrac{1}{2}\|x^*-A(x-y)\|_*^2\nonumber\\
&= \inf_{y\in X}q^*_{A_+}\big(x^*-\tfrac{1}{2}A(x-y)+\tfrac{1}{2}A^*(x+y)
\big)
+q_{A_+}(x-y)
+\tfrac{1}{2}\|y\|^2+\tfrac{1}{2}\|x^*-A(x-y)\|_*^2\nonumber\\
&=\inf_{y\in X} q^*_{A_+}\big(x^*-Ax+A_+(x+y)\big)
+q_{A_+}(x-y)
+\tfrac{1}{2}\|y\|^2+\tfrac{1}{2}\|x^*-A(x-y)\|_*^2\nonumber\\
&=\inf_{y \in X}q^*_{A_+}(x^*-Ax)+\scal{x+y}{x^*-Ax}+q_{A_+}(x+y)
+q_{A_+}(x-y)\nonumber\\ 
&\qquad\quad +\tfrac{1}{2}\|y\|^2+\tfrac{1}{2}\|x^*-A(x-y)\|_*^2\nonumber\\
&=\inf_{y \in X}q^*_{A_+}(x^*-Ax)+\scal{x+y}{x^*-Ax}+2q_{A_+}(x)
+2q_{A_+}(y)\nonumber\\ 
&\qquad\quad +\tfrac{1}{2}\|-\negthinspace
y\|^2+\tfrac{1}{2}\|x^*-A(x-y)\|_*^2\nonumber\nonumber\\
&\geq q^*_{A_+}(x^*-Ax) + \scal{x}{x^*}+ \inf_{y\in X}  \scal{y}{x^*-Ax} +
2q_{A_+}(y) + \scal{-y}{x^*-A(x-y)}\nonumber\\
&=q^*_{A_+}(x^*-Ax) + \scal{x}{x^*}+ \inf_{y\in X} 2q_{A_+}(y) + \scal{-y}{Ay}\nonumber\\
&= q^*_{A_+}(x^*-Ax) + \scal{x}{x^*}\nonumber\\
&= \mathcal{A}_A(x,x^*).
\end{align} 
}Hence $\mathcal{B}_A \geq \mathcal{A}_A$.
On the other hand, both $\mathcal{A}_A$ and $\mathcal{B}_A$
are autoconjugate (see Fact~\ref{GFF:2} and Fact~\ref{GFF:1}).
Altogether, Fact~\ref{f:auto}\ref{Au:4} implies Claim~3.

Finally, observe that Claims 1--3 yield the result. 
\end{proof}

\begin{example}
Suppose that $X$ is the Euclidean plane $\RR^2$, 
let $\theta \in\left[0,\tfrac{\pi}{2}\right[$, and set 
\begin{align}
A&=\begin{pmatrix}
\cos\theta & -\sin\theta \\
\sin \theta  & \cos\theta
\end{pmatrix} 
\quad\text{and}\quad
A_{\pi/2}=\begin{pmatrix}
0 & -1 \\
1  & 0
\end{pmatrix}.
\end{align} 
Then for every $(x,x^*)\in\RR^2\times\RR^2$, 
\begin{align}
\mathcal{A}_A(x,x^*)&=\mathcal{B}_A(x,x^*) =\mathcal{C}_A(x,x^*)\notag\\
&= \frac{1}{2\cos\theta}\|x^*-Ax\|^2+
\scal{x}{x^*}\notag\\
&=\frac{1}{2\cos\theta}\|x^*-(\sin\theta)
A_{\pi/2}x\|^2+\frac{\cos\theta}{2}\|x\|^2.
\end{align}
\end{example}
\begin{proof}
This follows from Theorem~\ref{t:main} since
$A_+ = (\cos\theta)\Id$, $q_{A_+} = (\cos\theta)\thalb\|\cdot\|^2$,
and $A_{\mathlarger\circ} = (\sin\theta)A_{\pi/2}$. 
\end{proof}

\section{Observations on Nonuniqueness}

\label{s:uniq}

Theorem~\ref{t:main} might nurture the conjecture that for
continuous linear monotone operators, all autoconjugate representers
coincide. This conjecture is false --- we shall provide a whole
family of distinct autoconjugate representers for the identity on $\RR$.
Our constructions rests on the following result.

\begin{proposition} \label{p:rrsp} 
Let $g\colon\RR\to\RX$ be such that 
\begin{equation} \label{e:rrsp}
(\forall x\in\RR)\quad g^*(-x)=g(x)\geq 0.
\end{equation}
Then
\begin{equation}\label{e:g0=0}
g(0)=0.
\end{equation}
Moreover, each of the following functions satisfies \eqref{e:rrsp}:
\begin{enumerate}
\item the indicator function $\displaystyle \iota_{\RP}\colon x\mapsto 
\begin{cases}
0, &\text{if $x\geq 0$;}\\
\pinf, &\text{if $x<0$;}
\end{cases} $
\item the energy function $\thalb|\cdot|^2$;
\item for $p>1$ and $q>1$ such that $\tfrac{1}{p}+\tfrac{1}{q}=1$,
the function 
$\displaystyle
x \mapsto
\begin{cases}
\tfrac{1}{p} x^p, &\text{if $x\geq 0$;}\\
\tfrac{1}{q} (-x)^q, &\text{if $x<0$.}
\end{cases}
$
\end{enumerate}
\end{proposition}
\begin{proof}
On the one hand, $g(0)\geq 0$.
On the other hand, 
$g(0)=g^*(-0) = g^*(0)=\sup_{y\in\RR} -g(y) = -\inf_{y\in\RR} g(y) \leq 0$.
Altogether, $g(0)=0$ and so \eqref{e:g0=0}.
It is straight-forward to verify that 
each of the given functions satisfies \eqref{e:rrsp}. 
\end{proof}

\begin{theorem} \label{t:rrsp}
Let $g\colon\RR\to\RX$ be such that for every $x\in\RR$, 
$g^*(-x)=g(x)\geq 0$, 
and set $q\colon \RR\to\RR\colon x\mapsto \thalb|x|^2$.
Then 
\begin{equation}
F\colon \RR^2 \to \RX\colon (x,y)\mapsto q\Big(\frac{x+y}{\sqrt{2}}\Big)
+ g\Big(\frac{x-y}{\sqrt{2}}\Big).
\end{equation}
is an autoconjugate representer for $\Id\colon \RR\to\RR\colon x\mapsto x$.
\end{theorem}
\begin{proof}
Let $(x,y)\in\RR^2$. Using the fact that $q^*=q$ and the assumption on $g$,
we see that 
\begin{align}
F^*(y,x)&=\sup_{(u,v)\in\RR^2} uy+ vx
-q\big(\tfrac{u+v}{\sqrt{2}}\big)-g\big(\tfrac{u-v}{\sqrt{2}}\big) \nonumber\\
&= \sup_{(u,v)\in\RR^2}\tfrac{u+v}{2}(x+y)-\tfrac{u-v}{2}(x-y)-
q\big(\tfrac{u+v}{\sqrt{2}}\big)-
g\big( \tfrac{u-v}{\sqrt{2}}\big) \nonumber\\
&=\sup_{(u,v)\in\RR^2} 
\tfrac{u+v}{\sqrt{2}}\tfrac{x+y}{\sqrt{2}}- 
\tfrac{u-v}{\sqrt{2}}\tfrac{x-y}{\sqrt{2}}-
q\big(\tfrac{u+v}{\sqrt{2}}\big)- 
g\big(\tfrac{u-v}{\sqrt{2}}\big)\nonumber\\
&=q^*\big(\tfrac{x+y}{\sqrt{2}}\big) +
g^*\big(-\tfrac{x-y}{\sqrt{2}}\big)\nonumber\\
&=q\big(\tfrac{x+y}{\sqrt{2}}\big) +
g\big(\tfrac{x-y}{\sqrt{2}}\big)\nonumber\\
&=F(x,y). 
\end{align}
Hence $F$ is autoconjugate.
In view of \eqref{e:g0=0}, we have 
$(x,y)\in\gra\big(G(F)\big)$
$\Leftrightarrow$
$y\in G(F)x$ 
$\Leftrightarrow$
$F(x,y)=xy$
$\Leftrightarrow$
$q\big((x+y)/\sqrt{2}\,\big)+g\big((x-y)/\sqrt{2}\,\big) = xy$
$\Leftrightarrow$
$\tfrac{1}{4}(x+y)^2 + g\big((x-y)/\sqrt{2}\,\big) = xy$
$\Leftrightarrow$
$\tfrac{1}{4}(x-y)^2 + g\big((x-y)/\sqrt{2}\,\big) = 0$
$\Leftrightarrow$
$x-y=0$
$\Leftrightarrow$
$(x,y)\in\gra(\Id)$. 
\end{proof}

\begin{remark}
Consider Theorem~\ref{t:rrsp}.
If we set $g=q=\thalb|\cdot|^2$, 
then $F = q \oplus q=
q_{\Id} \oplus  q_{\Id}^* = \mathcal{C}_{\Id}$
by Corollary~\ref{c:G}.
Thus, this pleasingly symmetric choice of $g$ gives rise
to $\mathcal{A}_{\Id} = \mathcal{B}_{\Id} = \mathcal{C}_{\Id}$.
Proposition~\ref{p:rrsp} provides other choices of $g$ that lead
to different autoconjugate representers for $\Id$.
\end{remark}

Having settled the nonuniqueness of autoconjugate representers, 
it is natural to ask
``What makes the autoconjugate representers of Theorem~\ref{t:main}
special?''
The next result provides a complete answer for a large class of 
linear operators.

\begin{theorem}\label{t:hoe}
Let $A\colon X\to X^*$ be continuous, linear, monotone, and
symmetric.
Let $F\colon X\times X^*\to\RX$ be such that $\ran A$ is closed.
Then
\begin{equation}
F = \mathcal{C}_A 
\quad\Leftrightarrow\quad
\begin{cases}
\text{$F$ is autoconjugate},\\
F(0,0)=0,\\
\big(\forall(x,y)\in X\times X\big)\;\;F(x,Ay)=F(y,Ax).
\end{cases}
\end{equation}
\end{theorem}
\begin{proof}
``$\Rightarrow$'': By Corollary~\ref{c:G}, 
$F$ is autoconjugate and $F(0,0) = (q_A\oplus q_A^*)(0,0) = 0$.
Let $x$ and $y$ be in $X$. 
Using \eqref{Alf:1}, we have
$F(x,Ay) = (q_A\oplus q_A^*)(x,Ay) = q_A(x) + q_A^*(Ay) = q_A(x) + q_A(y)
= q_A(y) + q_A^*(Ax) = (q_A\oplus q_A^*)(y,Ax) = F(y,Ax)$.

``$\Leftarrow$'':
Let $(x,x^*)\in X\times X^*$. We proceed by verifying the next two claims.

\emph{Claim~1:} $x^*\notin\ran A$ $\Rightarrow$ $F(x,x^*)=\pinf$.\\
Assume that $x^*\notin \ran A$. The Separation Theorem yields
$z\in X$ such that 
\begin{equation} \label{e:saucer:a}
\scal{z}{x^*}>0
\end{equation}
and $\max\scal{z}{\ran A}=0$. Since $A$ is symmetric, we deduce that
$Az=0$. This implies
$(\forall\rho\in\RR)$ 
$F(\rho z,0) = F(\rho z,A0) = F\big(0,A(\rho z)\big) = F(0,0)=0$.
Thus
\begin{align} 
(\forall\rho\in\RR) \quad 
F(x,x^*) &= F(x,x^*) + F(\rho z,0) = F(x,x^*)+F^*(0,\rho z)\notag\\
&\geq \scal{x}{0} + \scal{\rho z}{x^*} = \rho\scal{z}{x^*}.
\label{e:saucer:c}
\end{align}
In view of \eqref{e:saucer:a}, we see that Claim~1 
follows by letting $\rho\to\pinf$ in \eqref{e:saucer:c}.

\emph{Claim~2:} $x^*\in\ran A$ $\Rightarrow$ $F(x,x^*)\geq
\mathcal{C}_A(x,x^*)$.\\
Assume that $x^*\in\ran A$, say $x^*=Ay$. Then
$2F(x,x^*) = 2F(x,Ay) = F(x,Ay) + F(y,Ax) = F(x,Ay) + F^*(Ax,y)
\geq \scal{x}{Ax} + \scal{y}{Ay}$ and hence,
using \eqref{Alf:1},
\begin{equation}
F(x,x^*) \geq q_A(x) + q_A(y) = q_A(x) + q_A^*(Ay) = (q_A\oplus
q_A^*)(x,x^*).
\end{equation}
This and \eqref{e:symmG} yield Claim~2.

Note that Claim~1 and Claim~2 yield $F\geq \mathcal{C}_A$.
Therefore, Fact~\ref{f:auto}\ref{Au:4} implies that $F=\mathcal{C}_A$. 
\end{proof}

\section{Autoconjugate Representers for  $\partial(-\ln)$}

\label{s:-ln}

Theorem~\ref{t:main} showed that three ostensibly different
autoconjugate representers are in fact identical for continuous
linear monotone operators. It is tempting to consider a subdifferential
operator $\partial f$, and to compare $\mathcal{A}_{\partial f}$, 
$\mathcal{B}_{\partial f}$, and $f\oplus f^*$.
It turns out that these autoconjugate representers for $\partial f$
may all be different.
To aid in the construction of this example, it will be convenient 
to work in this section with the 
negative natural logarithm function
\begin{equation} \label{e:-ln}
f\colon \RR\to\RX \colon x\mapsto 
\begin{cases}
-\ln(x), &\text{if $x>0$;}\\
\pinf, &\text{if $x\leq 0$,}
\end{cases}
\end{equation}
and with the set
\begin{equation} \label{e:C}
C = \menge{(x,x^*)\in\RR\times\RR}{x^* \leq -\tfrac{1}{x}<0}.
\end{equation}
It is well known that
\begin{equation} \label{e:-ln*} 
(\forall x\in\RR)\quad
f^*(x) = -1+f(-x)
\end{equation}
and straight-forward to verify that 
\begin{align}
\tfrac{1}{\sqrt{2}}C 
&= \menge{(x,x^*)\in\RR\times\RR}{x^* \leq
-\tfrac{1}{2x}<0},\label{e:straight:1}\\
\tfrac{1}{2}C 
&= \menge{(x,x^*)\in\RR\times\RR}{x^* \leq -\tfrac{1}{4x}<0}
\label{e:straight:2},
\end{align}
and
\begin{equation} \label{e:straight:3} 
\tfrac{1}{\sqrt{2}}C 
\subsetneq
\tfrac{1}{{2}}C 
\subsetneq
\left]0,\pinf\right[ \times \left]\minf,0\right[.
\end{equation}
Furthermore,  \cite[Example~3.4]{BM} yields
\begin{equation} \label{e:montag:a}
\big(\forall (x,x^*)\in\RR\times\RR\big)\quad
F_{\partial f}(x,x^*) = 
\begin{cases}
1 - 2\sqrt{-xx^*}, &\text{if $x\geq 0$ and $x^*\leq 0$;}\\
\pinf, &\text{otherwise,}
\end{cases}
\end{equation}
and 
\begin{equation} \label{e:montag:b}
F_{\partial f}^{*\intercal} = -1+\iota_C.
\end{equation}

\begin{theorem} \label{t:different}
The functions
$\mathcal{A}_{\partial f}$, 
$\mathcal{B}_{\partial f}$, and
$f\oplus f^*$ 
have domains $\tfrac{1}{\sqrt{2}}C$, $\tfrac{1}{2}C$,
and $\left]0,\pinf\right[ \times \left]\minf,0\right[$, respectively.
Consequently, they are three different autoconjugate representers for
$\partial f$. 
\end{theorem}
\begin{proof}
Using \eqref{e:18}, \eqref{e:montag:a}, \eqref{e:montag:b} and
\eqref{e:straight:1}, we see that
\begin{align}
\dom \mathcal{A}_{\partial f} &= \menge{(x,\thalb x_1^* + \thalb x_2^*)
\in\RR\times\RR}{(x,x_1^*)\in\dom F_{\partial f}\;\text{and}\;
(x,x_2^*)\in\dom F_{\partial f}^{*\intercal}}\notag\\
&= \menge{(x,\thalb x_1^* + \thalb x_2^*)
\in\RR\times\RR}{x\geq 0,\; x_1^*\leq 0, \;\text{and}\;
(x,x_2^*)\in C}\notag\\
&= \menge{(x,x^*)\in\RR\times\RR}{x^*\leq -\tfrac{1}{2x}<0}\notag\\
&= \tfrac{1}{\sqrt{2}}C,
\end{align}
as claimed. Similarly, by \eqref{e:GFF:1}, \eqref{e:montag:a}, and
\eqref{e:montag:b}, 
\begin{align}
\dom\mathcal{B}_{\partial f} &= \thalb \dom F_{\partial f} + \thalb \dom
F_{\partial f}^{*\intercal}\notag\\
&= \thalb\big(\left[0,\pinf\right[\times \left]\minf,0\right]\big) + \thalb C\notag\\
&= \thalb C.
\end{align}
Furthermore, by \eqref{e:-ln} and \eqref{e:-ln*}, 
\begin{equation}
\dom (f\oplus f^*) = (\dom f)\times (\dom f^*) = \left]0,\pinf\right[
\times \left]\minf,0\right[.
\end{equation}
We thus have verified the statements concerning the domains.
Fact~\ref{GFF:2}, Fact~\ref{GFF:1}, and Example~\ref{ex:Gh} imply
that all three functions are autoconjugate representers for $\partial
f$. In view of \eqref{e:straight:3}, these functions are all different
since their domains are also all different. 
\end{proof}

\begin{remark}
Using \eqref{e:montag:a} and \eqref{e:montag:b}, one may verify
that
\begin{equation}
\big(\forall (x,x^*)\in\RR\times\RR\big)\quad
\mathcal{A}_{\partial f}(x,x^*) = 
\begin{cases}
-\sqrt{-1-2xx^*}, &\text{if $(x,x^*)\in \tfrac{1}{\sqrt{2}}C$;}\\
\pinf, &\text{otherwise.}
\end{cases}
\end{equation}
However, we do not have an explicit formula for $\mathcal{B}_{\partial f}$. 
\end{remark}

\section{Discontinuous Symmetric Operators}

\label{s:disc}

In this final section, we investigate discontinuous
symmetric operators. Specifically, we assume throughout this section that
$A\colon X\To X^*$ is maximal monotone, at most single-valued, $\dom A$ is
a linear subspace, and $A|_{\dom A}$ is linear and symmetric.
Put differently, we assume that 
\begin{equation}
A \colon \dom A \to X^* \quad\text{is linear, symmetric, and 
maximal montone.}
\end{equation}
It is convenient to extend the definition of $q_A$ in \eqref{e:diequad} 
to this more general setting via 
\begin{equation}
q_A \colon X\to\RR \colon
x \mapsto
\begin{cases}
\thalb\scal{x}{Ax}, &\text{if $x\in \dom A$;}\\
\pinf, &\text{otherwise.}
\end{cases}
\end{equation}
A key tool is the function 
\begin{equation}
f\colon X\to\RX \colon
x \mapsto \sup_{y\in \dom A}\scal{x}{Ay}-\thalb\scal{y}{Ay},
\end{equation}
which was introduced by Phelps and Simons.

\begin{fact}[Phelps-Simons]  \label{f:PheSim}
{\rm (See \cite{PheSim}.)}
The following hold.
\begin{enumerate}
\item \label{f:PheSim:lsc}
$f$ is proper, lower semicontinuous, and convex.
\item \label{f:PheSim:par}
$A = \partial f$. 
\item \label{f:PheSim:quad} 
$\dom A\subseteq \dom f \subseteq \cldom A$ and 
$(\forall x\in \dom A)$ $f(x) = \thalb\scal{x}{Ax}$.
\item \label{f:PheSim:dom}
$A$ is continuous $\Leftrightarrow$ $\dom A= X$ $\Leftrightarrow$
$\dom f= X$.
\end{enumerate}
\end{fact}

\begin{corollary} \label{c:PheSim}
The following hold.
\begin{enumerate}
\item \label{zahn:i} $f+\iota_{\dom A} = q_A$.
\item \label{zahn:iii} $f = q_A$ $\Leftrightarrow$ $\dom f = \dom A$.
\item \label{zahn:ii} $q_{A}^{**}  = f$.
\item \label{zahn:iv} If $A$ is one-to-one, then $f = q_{A^{-1}}^*$. 
\end{enumerate}
\end{corollary}
\begin{proof}
\ref{zahn:i}: Clear from Fact~\ref{f:PheSim}\ref{f:PheSim:quad}.
\ref{zahn:iii}: 
Since $\dom q_A = \dom A$, this item is a consequence of \ref{zahn:i}.
\ref{zahn:ii}: 
Using Fact~\ref{f:PheSim}\ref{f:PheSim:lsc}\&\ref{f:PheSim:par}
and 
a result by J.\ Borwein (see \cite[Theorem~1]{Borwein82} or 
\cite[Theorem~3.1.4(i)]{Zalinescu}), we see that
$f = f^{**} = (f+\iota_{\dom\partial f})^{**} = 
(f+\iota_{\dom A})^{**} = q_A^{**}$. 
\ref{zahn:iv}: If $A$ is one-to-one, then
\begin{equation}
(\forall x\in X)\quad
f(x) = \sup_{y^*\in\dom A^{-1}} \scal{x}{y^*}-\thalb\scal{A^{-1}y^*}{y^*}
= \sup_{y^*\in\dom q_{A^{-1}}} \scal{x}{y^*} - q_{A^{-1}}(y^*)
= q_{A^{-1}}^*(x).
\end{equation}
This completes the proof. 
\end{proof}

\begin{proposition} \label{p:domchar}
We have: $\dom A = \dom f$ $\Leftrightarrow$
every sequence $(x_n)_\nnn$ in $\dom A$
such that $(x_n)_\nnn$ and $(\scal{x_n}{Ax_n})_\nnn$ are convergent must
satisfy $\lim x_n \in \dom A$.
\end{proposition}
\begin{proof}
``$\Rightarrow$'': 
Assume that $(x_n)_\nnn$ is a sequence in $\dom A$ such 
that $(x_n)_\nnn$ converges to $x\in X$ and
$(\scal{x_n}{Ax_n})_\nnn$ is also convergent.
Using Fact~\ref{f:PheSim}\ref{f:PheSim:lsc} and
Corollary~\ref{c:PheSim}\ref{zahn:i}, we have 
$x \in \dom f$ and thus $x\in\dom A$. 

``$\Leftarrow$'': 
In view of Fact~\ref{f:PheSim}\ref{f:PheSim:quad}, it suffices to
show that $\dom f \subseteq \dom A$.
To this end, let $x\in \dom f$.
In view of Corollary~\ref{c:PheSim}\ref{zahn:ii}, 
there exists a sequence $(x_n)_\nnn$ in $\dom q_A= \dom A$ such
$x_n\to x$ and $\thalb\scal{x_n}{Ax_n} = q_A(x_n)\to f(x)$. 
By assumption, $x\in\dom A$, as required. 
\end{proof}

\begin{theorem} \label{t:potato}
Let $B\colon X^*\to X$ be continuous, linear, symmetric, monotone, and
one-to-one. Suppose that $A = B^{-1}$. Then
\begin{equation} \label{e:potato:b}
\mathcal{A}_A = 
q_A\oplus q_B
= (q_B^* + \iota_{\dom A})\oplus q_B.
\end{equation}
and
\begin{equation} \label{e:potato:c}
\mathcal{B}_A = \mathcal{A}_A^{**} = q_B^*\oplus q_B
\end{equation}
are both representers for $A$. Furthermore, 
$\mathcal{A}_A = \mathcal{B}_A$
$\Leftrightarrow$
$\dom q_B^{*} = \dom A$.
\end{theorem}
\begin{proof}
Since $q_A = q_B^*+\iota_{\dom A}$ by
Corollary~\ref{c:PheSim}\ref{zahn:i}\&\ref{zahn:iv}, 
it suffices to verify the left equality in \eqref{e:potato:b}. 
Let $(x,x^*)\in X\times X^*$.
Using \eqref{e:18} and Fact~\ref{f:Fitz}\ref{f:Fitz:ndif:1}, 
we see that
\begin{align}
\mathcal{A}_A(x,x^*) &= 
\inf_{y^*}
\tfrac{1}{2}F_A(x,x^*+y^*)+\tfrac{1}{2}F_A^*(x^*-y^*,x)\notag\\
&= \inf_{y^*}
\tfrac{1}{2}F_B(x^*+y^*,x)+\tfrac{1}{2}F_B^*(x,x^*-y^*)\notag\\
&= \inf_{y^*}
\tfrac{1}{2}F_B(x^*+y^*,x)+\tfrac{1}{2}\big(\iota_{\gra B}(x^*-y^*,x) +
\scal{x^*-y^*}{B(x^*-y^*)}\big). \label{e:potato:a}
\end{align}
If $x\notin\ran B = \dom A$, then \eqref{e:potato:a} shows that
$\mathcal{A}_A(x,x^*)=\pinf$, as required.
So assume that $x\in\ran B=\dom A$. 
In view of \eqref{e:potato:a}, \eqref{e:es:4}, and \eqref{Alf:1}, 
we deduce that
\begin{align}
\mathcal{A}_A(x,x^*) 
&= \tfrac{1}{2}F_B(2x^*-Ax,x)+\tfrac{1}{2}\scal{x}{Ax}\notag\\
&= q_B^*\big(\thalb x + \thalb B(2x^*-Ax)\big)+q_A(x)\notag\\
&= q_B^*(Bx^*)+q_A(x)\notag\\
&= q_B(x^*)+q_A(x).
\end{align}
Hence \eqref{e:potato:b} holds.
Using \eqref{e:potato:b},
Corollary~\ref{c:PheSim}\ref{zahn:ii}\&\ref{zahn:iv}, \eqref{e:symmG}, and
Theorem~\ref{t:main}, we see that 
\begin{equation}
\mathcal{A}_A^{**} = (q_A\oplus q_B)^{**} 
= q_{A}^{**} \oplus q_{B}^{**} = q_{B}^*\oplus q_B = 
(q_B\oplus q_{B}^*)^* = \mathcal{B}_B^* = \mathcal{B}_B^\intercal 
=\mathcal{B}_{B^{-1}} = \mathcal{B}_A, 
\end{equation}
so that \eqref{e:potato:c} holds. Furthermore, 
$\mathcal{A}_A = \mathcal{B}_A$ $\Leftrightarrow$
$q_A \oplus q_B = q_{A}^{**}\oplus q_B$ $\Leftrightarrow$
$q_A = q_{A}^{**}$ $\Leftrightarrow$ $\dom q_{A}^{**} = \dom A$
$\Leftrightarrow$ $\dom q_B^* = \dom A$ by Corollary~\ref{c:PheSim}. 
\end{proof}

\begin{example} \label{ex:potato}
Suppose that $X$ is the Hilbert space $\ell_2(\NN)$ of square-summable
sequences; thus, $X^*=X$.
Set 
\begin{equation}
B\colon X\to X \colon (\xi_k)_{k\in\NN} \mapsto
(\tfrac{1}{k}\xi_k)_{k\in\NN}
\end{equation}
and suppose that $A= B^{-1}$. 
Then $\ran B = \dom A$ is dense in $X$, but it is 
not closed (since, e.g., 
$(\tfrac{1}{k})_{k\in\NN}\in X\smallsetminus(\ran B)$). 
Now set
\begin{equation}
x = \big(\tfrac{1}{k^{4/3}}\big)_{k\in\NN}
\quad\text{and}\quad
(\forall \nnn)\;\;
x_n =
\big(\tfrac{1}{1^{4/3}},\tfrac{1}{2^{4/3}},\ldots,\tfrac{1}{n^{4/3}},0,0,\ldots\big).
\end{equation}
On the one hand, $(x_n)_\nnn$ lies in $\dom A$ and 
$x_n\to x\in X\smallsetminus(\dom A)$. 
On the other hand, $\scal{x_n}{Ax_n} = \sum_{k=1}^{n}
\tfrac{1}{k^{4/3}}\tfrac{k}{k^{4/3}} =
\sum_{k=1}^{n} \tfrac{1}{k^{5/3}} \to \zeta(5/3)\in\RR$. 
Altogether, Proposition~\ref{p:domchar} implies that 
$\dom A\subsetneq \dom q_B^*$. Therefore, by Theorem~\ref{t:potato},
$\mathcal{A}_A$ is neither lower semicontinuous nor equal to
$\mathcal{B}_A$. While $\mathcal{A}_A$ is still a representer for $A$,
it cannot be autoconjugate. 
\end{example}

\begin{remark} Several comments are in order. 
\begin{enumerate}
\item Without the constraint qualification, Fact~\ref{GFF:2} fails 
(see Example~\ref{ex:potato}, where $\dom A$ is a subspace that is not
closed). 
\item 
It is conceivable that $\mathcal{A}_A^{**}$ is always an
autoconjugate representer for $A$ --- 
this would sharpen Fact~\ref{GFF:2} and it would be consistent with
Theorem~\ref{t:potato}. 
\item \label{r:potato:whether} 
Suppose that $B$ is as in Theorem~\ref{t:potato}, that $A=B^{-1}$, 
and
that $\dom A = \ran B$ is a dense subspace of $X$ with $\dom A\neq X$.

We do not know whether $(\dom f)\smallsetminus (\dom A)\neq\emp$ must
hold (as it does in Example~\ref{ex:potato}), i.e.\ (see
Proposition~\ref{p:domchar}),
whether there exists a sequence $(x_n)_\nnn$ in $\dom A$ such that
$(x_n)_\nnn$  converges to some point $x\in X\smallsetminus(\dom A)$, yet
$(\scal{x_n}{Ax_n})_\nnn$ converges to a real number. 

In contrast, there does exist a point $x\in X\smallsetminus(\dom A)$ such that
every sequence $(x_n)_\nnn$ in $\dom A$ converging to $x$ must have 
$\scal{x_n}{Ax_n} \to \pinf$.
(Indeed, since $\dom A\neq X$, it follows from
Fact~\ref{f:PheSim}\ref{f:PheSim:dom} that $\dom f\neq X$. Take $x\in
X\smallsetminus(\dom f)$ and assume that $(x_n)_\nnn$ lies in $\dom A$ and
converges to $x$. Then $\pinf = f(x)\leq\varliminf f(x_n) = \varliminf
\thalb \scal{x_n}{Ax_n}$ by
Fact~\ref{f:PheSim}\ref{f:PheSim:lsc}\&\ref{f:PheSim:quad}). 
Thus for 
every sequence $(x_n^*)_\nnn$ in $X^*$ such that 
$Bx_n^*\to x\notin\ran B$, it follows that 
$\|Bx_n^*\|\cdot\|x_n^*\|_* \geq \scal{Bx_n^*}{x_n^*} = 
\scal{Bx_n^*}{A(Bx_n^*)}\to\pinf$. 
Since $0\in\dom f$ and so $x\neq 0$, we deduce 
$\|x_n^*\|_* \to\pinf$, which is a 
well known result from Functional Analysis 
(see \cite[Corollary~17.G]{Holmes}). 
\end{enumerate}
\end{remark}

\section*{Acknowledgment}
Heinz Bauschke was partially supported by the Canada Research Chair Program
and by the Natural Sciences and
Engineering Research Council of Canada.  
Xianfu Wang was partially supported by the Natural
Sciences and Engineering Research Council of Canada.

\end{document}